\documentclass[12pt]{article}

\usepackage{mathrsfs}
\usepackage{amsfonts}
\usepackage{tikz}
\usepackage{amssymb,bm}
\usepackage{amsmath}
\usepackage{amsthm}
\usepackage{bbding}
\usetikzlibrary{matrix,shapes,snakes}

\allowdisplaybreaks
\newtheorem{thm}{Theorem}[section]
\newtheorem{defi}[thm]{Definition}
\newtheorem{lem}[thm]{Lemma}
\newtheorem{cor}[thm]{Corollary}
\newtheorem{prop}[thm]{Proposition}
\newtheorem{rmk}[thm]{Remark}
\newtheorem{eg}[thm]{Example}

\def\proof{\vskip 1mm\noindent{\it Proof.}\quad}
\newcommand{\qedd}{\hspace*{\fill}$\Box$\medskip}

\def\deg{\hbox{\rm{deg\,}}}
\def\rk{\hbox{\rm{rk\,}}}
\def\det{\hbox{\rm{det\,}}}
\def\gal{\hbox{\rm{Gal}}}
\def\im{\hbox{\rm{Im\,}}}
\renewcommand{\vec}[1]{\bm{#1}}
\newcommand{\DF}[2]{{\displaystyle\frac{#1}{#2}}}
\def\tr{\hbox{\rm{tr}}}
\def\Tr{\hbox{\rm{Tr}}}
\def\hom{\hbox{\rm{Hom}}}
\def\endm{\hbox{\rm{End}}}

\begin{document}

\title{Linearized polynomials over finite fields revisited\thanks{Partially supported
by National Basic Research Program of China (2011CB302400).}}

\author{Baofeng Wu\thanks{Key Laboratory of Mathematics
Mechanization, AMSS, Chinese Academy of Sciences,
 Beijing 100190,  China. Email: wubaofeng@amss.ac.cn}, Zhuojun Liu\thanks{Key Laboratory of Mathematics
Mechanization, AMSS, Chinese Academy of Sciences,
 Beijing 100190,  China. Email: zliu@mmrc.iss.ac.cn}}
 \date{}

\maketitle

\begin{abstract}
We give new characterizations of the algebra
$\mathscr{L}_n(\mathbb{F}_{q^n})$ formed by all linearized
polynomials over the finite field $\mathbb{F}_{q^n}$ after briefly
surveying some known ones. One isomorphism we construct is between
$\mathscr{L}_n(\mathbb{F}_{q^n})$ and the composition algebra
$\mathbb{F}_{q^n}^\vee\otimes_{\mathbb{F}_{q}}\mathbb{F}_{q^n}$. The
other isomorphism we construct is between
$\mathscr{L}_n(\mathbb{F}_{q^n})$ and the so-called Dickson matrix
algebra $\mathscr{D}_n(\mathbb{F}_{q^n})$. We also further study the
relations between a linearized polynomial and its associated Dickson
matrix, generalizing a well-known criterion of Dickson on linearized
permutation polynomials. Adjugate polynomial of a linearized
polynomial is then introduced, and connections between them are
discussed. Both of the new characterizations can bring us more
simple approaches to establish a special form of representations of
linearized polynomials proposed recently by several authors.
Structure of the subalgebra $\mathscr{L}_n(\mathbb{F}_{q^m})$ which
are formed by all linearized polynomials over a subfield
$\mathbb{F}_{q^m}$ of $\mathbb{F}_{q^n}$ where $m|n$ are also
described.\vskip .5em


\end{abstract}


\section{Introduction}\label{secintro}

Let $\mathbb{F}_q$ and $\mathbb{F}_{q^n}$ be the finite fields with
$q$ and $q^n$ elements respectively, where $q$ is a prime or a prime
power. Polynomials over $\mathbb{F}_{q^n}$ of the form
\begin{equation}\label{lp}
L(x)=\sum_{i=0}^{t}a_ix^{q^i},~t\in \mathbb{N}
\end{equation}
are often known as linearized polynomials. Such special kinds of
polynomials can induce linear transformations of $\mathbb{F}_{q^n}$
over $\mathbb{F}_q$. Considered as maps between finite fields,
linearized polynomials are always taken as
\begin{equation}\label{mlp}
L(x)=\sum_{i=0}^{n-1}a_ix^{q^i}\in \mathbb{F}_{q^n}[x]/(x^{q^n}-x).
\end{equation}

We denote by $\mathscr {L}(\mathbb{F}_{q^n})$ and $\mathscr
{L}_n(\mathbb{F}_{q^n})$  the set of all linearized polynomials in
the form (\ref{lp}) and (\ref{mlp}) respectively. Equipped with the
operations of addition and composition of polynomials in
$\mathbb{F}_{q^n}[x]$ and $\mathbb{F}_{q^n}[x]/(x^{q^n}-x)$
respectively, and scalar multiplication by elements in
$\mathbb{F}_q$, $\mathscr {L}(\mathbb{F}_{q^n})$ and $\mathscr
{L}_n(\mathbb{F}_{q^n})$ both
 form non-commutative $\mathbb{F}_q$-algebras. We
 mainly focus on the algebra $\mathscr
{L}_n(\mathbb{F}_{q^n})$ in this paper, and use ``$\circ$" to denote
the multiplication in it. Let $\mathscr
{L}_n(\mathbb{F}_{q^m})\subset\mathscr {L}_n(\mathbb{F}_{q^n})$ be
the set of linearized polynomials with coefficients in
$\mathbb{F}_{q^m}$, where $\mathbb{F}_{q^m}$ is the subfield of
$\mathbb{F}_{q^n}$ with $q^m$ elements (which implies $m|n$). It is
clear that $\mathscr {L}_n(\mathbb{F}_{q^m})$ is an
$\mathbb{F}_q$-subalgebra of $\mathscr {L}_n(\mathbb{F}_{q^n})$
under the operations mentioned above.

Linearized polynomials were firstly studied by Ore in \cite{ore33}
after his work on the theory of non-commutative polynomials
\cite{ore33i}. From \cite{ore33}, we can easily get the
$\mathbb{F}_q$-algebra isomorphism
\[\mathscr {L}_n(\mathbb{F}_{q^n})\cong\mathbb{F}_{q^n}[x; \sigma]/(x^n-1),\]
where $\sigma $ is the Frobenius automorphism of $\mathbb{F}_{q^n}$
over $\mathbb{F}_{q}$, i.e. $\sigma(x)=x^q$ for
$x\in\mathbb{F}_{q^n}$, and $\mathbb{F}_{q^n}[x; \sigma]$ is the
so-called skew-polynomial ring. Actually, it is clear that every
linear transformation of $\mathbb{F}_{q^n}$ over $\mathbb{F}_{q}$
can be induced by a linearized polynomial, thus we have
\[\mathscr {L}_n(\mathbb{F}_{q^n})\cong\mathscr {M}_n(\mathbb{F}_{q}),\]
where $\mathscr {M}_n(\mathbb{F}_{q})$ is the algebra formed by all
$n\times n$ matrices over $\mathbb{F}_{q}$. This isomorphism was
constructed explicitly by Carlitz in \cite{carl}. In addition,
$\mathscr {L}_n(\mathbb{F}_{q^n})$ is also proved to be isomorphic
to $\mathbb{F}_{q^n}[G]$, the so-called semi-linear group algebra,
where $G=\gal(\mathbb{F}_{q^n}/\mathbb{F}_{q})$. All these results
give characterizations of the algebra $\mathscr
{L}_n(\mathbb{F}_{q^n})$.

In addition to the known ones, we give two new approaches to
characterize the algebra $\mathscr {L}_n(\mathbb{F}_{q^n})$ in this
paper. One isomorphism we construct is
\begin{equation}\label{compiso}
\mathscr{L}_n(\mathbb{F}_{q^n})\cong
\mathbb{F}_{q^n}^\vee\otimes_{\mathbb{F}_{q}}\mathbb{F}_{q^n},
\end{equation}
where $\mathbb{F}_{q^n}^\vee$ is the dual space of
$\mathbb{F}_{q^n}$ over $\mathbb{F}_{q}$, and
$\mathbb{F}_{q^n}^\vee\otimes_{\mathbb{F}_{q}}\mathbb{F}_{q^n}$ is
the so-called composition algebra on $\mathbb{F}_{q^n}$. In fact,
this is a special case of the composition algebra mentioned in
\cite{greub}. The other isomorphism we construct is
\begin{equation}\label{dmatiso}
\mathscr{L}_n(\mathbb{F}_{q^n})\cong
\mathscr{D}_n(\mathbb{F}_{q^n}),
\end{equation}
where $\mathscr{D}_n(\mathbb{F}_{q^n})$ is an algebra formed by all
$n\times n$ matrices over $\mathbb{F}_{q^n}$ of the form
\begin{equation}
\label{dmat}\begin{pmatrix}
a_0&a_1&\dots&a_{n-1}\\
a_{n-1}^q&a_0^q&\dots&a_{n-2}^q\\
\vdots&\vdots&&\vdots\\
a_1^{q^{n-1}}&a_2^{q^{n-1}}&\dots&a_0^{q^{n-1}}
\end{pmatrix},
\end{equation}
which are called Dickson matrices.

A well known result of Dickson indicates that
$L(x)=\sum_{i=0}^{n-1}a_ix^{q^i}\in\mathscr{L}_n(\mathbb{F}_{q^n})$
is a linearized permutation polynomial if and only if the Dickson
matrix associated to it, i.e. the Dickson matrix with the
coefficients $a_0$, $a_1$, $\ldots$, $a_{n-1}$ as entries of the
first row, is non-singular \cite{lidl}. As a matter of fact, a
linearized permutation polynomial is just a linearized polynomial of
rank $n$. By the rank of a linearized polynomial, we mean the rank
of it as a linear transformation of $\mathbb{F}_{q^n}$ over
$\mathbb{F}_{q}$, i.e. the dimension of the image space over
$\mathbb{F}_{q}$. We can generalize Dickson's result to that the
rank of a linearized polynomial is equal to the rank of the Dickson
matrix associated to it. Furthermore, we find that the adjugate
matrix of a Dickson matrix is also a Dickson matrix, thus we can
define an adjugate polynomial for every linearized polynomial. From
this concept, we can derive more properties of linearized
polynomials of rank $n$ and $n-1$ respectively.

Recently, it was proved in \cite{lings} that for a fixed ordered
basis $\{\beta_i\}_{i=0}^{n-1}$,  any linearized polynomial
$L(x)\in\mathscr{L}_n(\mathbb{F}_{q^n})$ of rank $k$ can be
represented as
\[L(x)=\sum_{i=0}^{n-1}\tr (\beta_i x)\alpha_i\]
or \[L(x)=\sum_{i=0}^{n-1}\tr (\alpha'_i x)\beta_i,\]
 where $\{\alpha_i\}_{i=0}^{n-1}$ and $\{\alpha'_i\}_{i=0}^{n-1}$ are two order sets of elements
in $\mathbb{F}_{q^n}$ both of rank $k$ over $\mathbb{F}_{q}$, and
``$\tr$" is the trace function from $\mathbb{F}_{q^n}$ to
$\mathbb{F}_{q}$. We notice that that these kinds of representations
of linearized polynomials are just direct consequences of the
isomorphism (\ref{compiso}), thus we rediscover them via a more
simple approach. Besides, we can also establish them simply
according to the isomorphism (\ref{dmatiso}) and properties of
Dickson matrices.

As to the subalgebra $\mathscr {L}_n(\mathbb{F}_{q^m})$,  in
\cite{ore34} Ore derived that
\[\mathscr {L}_n(\mathbb{F}_{q})\cong\mathbb{F}_{q}[x]/(x^n-1),\]
and in \cite{brawley}, Brawley et al. completed the problem. They
derived the isomorphism
\begin{equation}\label{subalg}
\mathscr {L}_n(\mathbb{F}_{q^m})\cong\mathscr
{M}_m\big(\mathbb{F}_{q}[x]/(x^t-1)\big),
\end{equation}
where $n=mt$ and $\mathscr {M}_m\big(\mathbb{F}_{q}[x]/(x^t-1)\big)$
is the algebra formed by all $m\times m$ matrices over the residue
ring $\mathbb{F}_{q}[x]/(x^t-1)$. We can also give descriptions to
the structure of $\mathscr {L}_n(\mathbb{F}_{q^m})$ via constructing
isomorphisms between it and subalgebras of
$\mathbb{F}_{q^n}^\vee\otimes_{\mathbb{F}_{q}}\mathbb{F}_{q^n}$ and
$\mathscr{D}_n(\mathbb{F}_{q^n})$ respectively.

In this paper, we revisit linearized polynomials over finite fields
and get some new results related to them. The rest of the paper is
organized as follows. In Section 2, we briefly survey the known
characterizations of the algebra $\mathscr{L}_n(\mathbb{F}_{q^n})$.
In Section 3, we propose the composition algebra approach to
characterize $\mathscr{L}_n(\mathbb{F}_{q^n})$. In Section 4, we
propose the Dickson matrix algebra approach to characterize
$\mathscr{L}_n(\mathbb{F}_{q^n})$, and derive more relations between
linearized polynomials and their associated Dickson matrices. In
Section 5, the trace form representations of linearized polynomials
are rediscovered through more simple approaches and studied further.
In Section 6, descriptions to the subalgebra $\mathscr
{L}_n(\mathbb{F}_{q^m})$ are discussed. Concluding remarks are given
in Section 7.

Throughout this paper, we fix an ordered basis
$\{\beta_i\}_{i=0}^{n-1}$ of $\mathbb{F}_{q^n}$ over
$\mathbb{F}_{q}$, and denote its dual basis by
$\{\beta^*_i\}_{i=0}^{n-1}$, i.e.
$\tr(\beta_i\beta_j^*)=\delta_{ij}$, $0\leq i,~j\leq n-1$,  where
$\delta$ is the Kronecker delta. When $\{\beta_i\}_{i=0}^{n-1}$ is a
normal basis, we assume $\beta_i=\beta^{q^i}$, $0\leq i\leq n-1$,
where $\beta$ is the normal basis generator with dual basis
generator $\beta^*$. For a set of elements $\{\alpha_i\}_{i=0}^{s}$,
we denote by $\rk_{\mathbb{F}_{q}}\{\alpha_i\}_{i=0}^{s}$ to be its
rank over $\mathbb{F}_{q}$. The rank of a linearized polynomial
$L(x)$ and a matrix $M$ are denoted to be $\rk L$ and $\rk M$
respectively. For a matrix $\big(a_{ij}\big)_{0\leq i\leq
n-1,\,0\leq j\leq n-1}$ with the $(i,j)$-th entry $a_{ij}$, we
sometimes use $\big(a_{ij}\big)$ to represent it for simplicity. For
example, we always denote the Dickson matrix in (\ref{dmat}) by
$\left(a_{j-i}^{q^i}\right)$ (subscripts reduced modulo $n$).


\section{Known characterizations of $\mathscr
{L}_n(\mathbb{F}_{q^n})$}\label{secalg}

As mentioned in Section 1, the $\mathbb{F}_{q}$-algebra
$\mathscr{L}_n(\mathbb{F}_{q^n})$ can be characterized through
various approaches. For completeness of this paper, we firstly
recall the known ones briefly in this part.

\subsection{The skew-polynomial ring approach}\label{ssecskewpoly}

Let $\sigma$ be the Frobenius automorphism of $\mathbb{F}_{q^n}$
over $\mathbb{F}_{q}$ and $\mathbb{F}_{q^n}[x; \sigma]$ be the set
of all (skew-)polynomials over $\mathbb{F}_{q^n}$. Define the
addition of polynomials in the standard manner and multiplication of
polynomials by distribute law and
\[(ax^i)(bx^j)=a\sigma^i(b)x^{i+j}=ab^{q^i}x^{i+j},~~a,~b\in\mathbb{F}_{q^n}.\]
With these operations, $\mathbb{F}_{q^n}[x; \sigma]$ forms a ring
called the skew-polynomial ring over $\mathbb{F}_{q^n}$ with
automorphism $\sigma$. This is a special case of Ore's
non-commutative polynomial ring proposed in \cite{ore33i}. It can be
proved that $\mathbb{F}_{q^n}[x; \sigma]$ is a principle ideal
domain, i.e. a ring with no zero divisors and whose left and right
ideals are all principle. $\mathbb{F}_{q^n}[x; \sigma]$ also becomes
an $\mathbb{F}_{q}$-algebra under scalar multiplication of elements
in $\mathbb{F}_{q}$.

Define
\begin{eqnarray*}
\phi: ~\mathbb{F}_{q^n}[x; \sigma]&\longrightarrow&
\mathscr{L}(\mathbb{F}_{q^n})\\
\sum_{i=0}^{t}a_ix^{i}&\longmapsto& \sum_{i=0}^{t}a_ix^{q^i}
\end{eqnarray*}
Note that $(ax^i)(bx^j)=ab^{q^i}x^{i+j}$ is mapped to
$ab^{q^i}x^{q^{i+j}}=(ax^{q^i})\circ(bx^{q^j})$. The following
theorem is straightforward.

\begin{thm}
The above map $\phi$ defines an algebra isomorphism
\[\mathscr{L}(\mathbb{F}_{q^n})\cong\mathbb{F}_{q^n}[x; \sigma].\]
\end{thm}

It is easy to find that $(x^n-1)$ and $(x^{q^n}-x)$ are two-sided
ideals of $\mathbb{F}_{q^n}[x; \sigma]$ and
$\mathscr{L}(\mathbb{F}_{q^n})$ respectively, and
$\phi(x^n-1)=x^{q^n}-x$, so a map $\bar{\phi}$ between quotient
algebras can be induced by $\phi$. Hence we have:
\begin{thm}[\cite{ore33}]
\[\mathscr{L}_n(\mathbb{F}_{q^n})=\mathscr{L}(\mathbb{F}_{q^n})/(x^{q^n}-x)\cong\mathbb{F}_{q^n}[x; \sigma]/(x^n-1).\]
\end{thm}

\begin{rmk}\label{skewgcrd}
Since $\mathbb{F}_{q^n}[x; \sigma]$ is a right Euclidean domain, the
greatest common right divisor (gcrd) of two skew-polynomials can be
computed \cite{ore33i}. It is easy to prove that $\rk L=n-\deg
\mbox{\rm{gcrd}}(\phi^{-1}(L(x)), x^n-1)$ for any
$L(x)\in\mathscr{L}_n(\mathbb{F}_{q^n})$ (viewed as an element of
$\mathscr{L}(\mathbb{F}_{q^n})$ formally), where the degree of a
skew-polynomial $f(x)$, $\deg f$, is defined to be the degree of it
as a usual polynomial. This supplies an approach to get rank of a
given linearized polynomial. Algorithms for computing gcrd of two
skew-polynomials and their complexity can be found in \cite{gies}.
As a special case, we know that $\rk L=n-1$ if and only if $\deg
\mbox{\rm{gcrd}}(\phi^{-1}(L(x)), x^n-1)=1$. This is equivalent to
say $L(x)$ is of the form
\[L(x)=L_1(x)\circ(x^q-ax),\]
where $a\in\big(\mathbb{F}_{q^n}\big)^{q-1}$ and $\deg L_1=(\deg
L)/q$ (the degree of a linearized polynomial in
$\mathscr{L}_n(\mathbb{F}_{q^n})$ is defined in the usual sense
viewing it as a polynomial in $\mathscr{L}(\mathbb{F}_{q^n})$
formally). This will supply an answer to the open problem proposed
in \cite{charp} to an extent. We will discuss this in detail in
Section \ref{seccompalg}.
\end{rmk}

\subsection{The semi-linear group algebra approach}\label{ssecgpring}

Let $G=\gal(\mathbb{F}_{q^n}/\mathbb{F}_{q})=\langle\sigma\rangle$
be the Galois group of  $\mathbb{F}_{q^n}$ over $\mathbb{F}_{q}$ and
$\mathbb{F}_{q^n}[G]$ be the vector space generated by elements in
$G$ over $\mathbb{F}_{q^n}$. Define the multiplication operation in
$\mathbb{F}_{q^n}[G]$ by the distribute law and
\[(a\sigma^i)(b\sigma^j)=a\sigma^i(b)\sigma^{i+j}=ab^{q^i}\sigma^{i+j},~~a,~b\in\mathbb{F}_{q^n}.\]
With this operation, $\mathbb{F}_{q^n}[G]$ forms a non-commutative
ring called the semi-linear group ring over $\mathbb{F}_{q^n}$. It
is also an $\mathbb{F}_{q}$-algebra under scalar multiplication by
elements in $\mathbb{F}_{q}$.

By the definition of the multiplication in $\mathbb{F}_{q^n}[G]$, it
is easy to see that the map
\begin{eqnarray*}
\mathbb{F}_{q^n}[G]&\longrightarrow&
\mathscr{L}_n(\mathbb{F}_{q^n})\\
\sum_{i=0}^{n-1}a_i\sigma^{i}&\longmapsto&
\sum_{i=0}^{n-1}a_ix^{q^i}
\end{eqnarray*}
is an algebra homomorphism. By comparing dimensions, we can get
\begin{thm}[\cite{mcdonald}]
\[\mathscr{L}_n(\mathbb{F}_{q^n})\cong\mathbb{F}_{q^n}[G].\]
\end{thm}

\subsection{The matrix algebra approach}\label{ssecmatalg}

Let $\mathscr {M}_n(\mathbb{F}_{q})$ be the matrix algebra over
$\mathbb{F}_{q}$ and $\endm(\mathbb{F}_{q^n})$ be the algebra of
linear transformations of $\mathbb{F}_{q^n}$ over $\mathbb{F}_{q}$.
It is well known that every linear transformation of
$\mathbb{F}_{q^n}$ over $\mathbb{F}_{q}$ can be induced by a
linearized polynomial over $\mathbb{F}_{q^n}$. This is because
$\endm(\mathbb{F}_{q^n})$ is an $\mathbb{F}_{q}$-vector space, which
becomes an $\mathbb{F}_{q^n}$-vector space under the scalar
multiplication
\[(aT)(x)=aT(x), ~\forall x\in \mathbb{F}_{q^n}\]
for any $a\in \mathbb{F}_{q^n}$ and $T\in\endm(\mathbb{F}_{q^n})$.
Note that $\{\sigma^i\in
\endm(\mathbb{F}_{q^n})|\;0\leq i\leq n-1\}$ is linearly
independent over $\mathbb{F}_{q^n}$ since if not so, the polynomial
$\sum_{i=0}^{n-1}a_i\sigma^{i}(x)=\sum_{i=0}^{n-1}a_ix^{q^i}$ will
admit $q^n$ roots for some $a_0,\cdots, a_{n-1}\in\mathbb{F}_{q^n}$
which are not all zero. As
$\dim_{\mathbb{F}_{q}}\endm(\mathbb{F}_{q^n})=n^2$ and
$\dim_{\mathbb{F}_{q}}\mathbb{F}_{q^n}=n$, we have
\[\dim_{\mathbb{F}_{q^n}}\endm(\mathbb{F}_{q^n})=n,\]
which implies that for any $T\in\endm(\mathbb{F}_{q^n})$, there
exist $a_0,\cdots, a_{n-1}\in\mathbb{F}_{q^n}$ such that
$T=\sum_{i=0}^{n-1}a_i\sigma^{i}$. Define a map
\begin{eqnarray*}
\endm(\mathbb{F}_{q^n})&\longrightarrow&
\mathscr{L}_n(\mathbb{F}_{q^n})\\
\sum_{i=0}^{n-1}a_i\sigma^{i}&\longmapsto&
\sum_{i=0}^{n-1}a_ix^{q^i}.
\end{eqnarray*}
It is easy to find that this map is bijective.

Moreover, it is straightforward to verify that the map defined above
is an $\mathbb{F}_{q}$-algebra isomorphism. From linear algebra we
finally
 get:
\begin{thm}[\cite{mcdonald}]\label{isomatalg}
  \[\mathscr{L}_n(\mathbb{F}_{q^n})\cong \mathscr{M}_n(\mathbb{F}_{q}).\]
\end{thm}

In \cite{carl} Carlitz explicitly constructed the matrix
corresponding to a linearized permutation polynomial by fixing a
normal basis of $\mathbb{F}_{q^n}$ over $\mathbb{F}_{q}$. However,
we find that the matrix corresponding to any linearized polynomial
can be explicitly constructed under any fixed basis with the aid of
the Dickson matrix associated to the given polynomial. We will talk
about this in detail in Section \ref{secdmatalg}.


\section{Characterize $\mathscr{L}_n(\mathbb{F}_{q^n})$ via the composition algebra}\label{seccompalg}

Let $E$ be a vector space over a field $\mathbb{K}$, and $E^\vee$ be
its dual space, i.e. $E^\vee=\hom_\mathbb{K}(E,\mathbb{K})$. Define
a multiplication in the tensor space $E^\vee\otimes_\mathbb{K} E$ by
setting
\[(l_1\otimes x_1)(l_2\otimes x_2)=l_1(x_2)l_2\otimes x_1\]
for any $l_1,~l_2\in E^\vee$ and $x_1,~x_2\in E$, and expanding it
by linearity to the whole space. It is easy to verify that this
multiplication makes $E^\vee\otimes_\mathbb{K} E$ into an associate
(non-commutative) $\mathbb{K}$-algebra called the composition
algebra \cite{greub}.

Consider the map $\Lambda: ~E^\vee\otimes_\mathbb{K}
E\longrightarrow
\endm_\mathbb{K}(E)$ defined by $\Lambda(l\otimes x)(y)=l(y)x$ for
any $y\in E$. It can be verified that $\Lambda$ is an injective
algebra homomorphism. Thus $\Lambda$ is an algebra isomorphism when
$E$ is furthermore of a finite dimension.

Now we let $E=\mathbb{F}_{q^n}$, which is an $n$-dimensional
$\mathbb{F}_{q}$-vector space, and consider the composition algebra
$\mathbb{F}_{q^n}^\vee\otimes_{\mathbb{F}_{q}}\mathbb{F}_{q^n}$.
Recall the isomorphism between $\mathscr{L}_n(\mathbb{F}_{q^n})$ and
$\endm(\mathbb{F}_{q^n})$, we get the following theorem
straightforwardly from the above discussions.
\begin{thm}\label{compo}
\[\mathscr{L}_n(\mathbb{F}_{q^n})\cong
\mathbb{F}_{q^n}^\vee\otimes_{\mathbb{F}_{q}}\mathbb{F}_{q^n}.\]
\end{thm}

In fact, the isomorphism in Theorem \ref{compo} can be constructed
directly. Let $T_\alpha\in\mathbb{F}_{q^n}^\vee$ be the linear
function defined by
\[T_\alpha(x)=\tr(\alpha x),~~\forall x\in \mathbb{F}_{q^n},\]
for any $\alpha\in\mathbb{F}_{q^n}$. As is well known that
\[\mathbb{F}_{q^n}^\vee=\{T_\alpha|\,\alpha\in\mathbb{F}_{q^n}\}\]
(a short proof is like this: define a scalar multiplication by
elements in $\mathbb{F}_{q^n}$ as $(\alpha T)(x)=T(\alpha x)$,
$\forall x\in\mathbb{F}_{q^n}$, for any $T\in\mathbb{F}_{q^n}^\vee$,
which makes $\mathbb{F}_{q^n}^\vee$ into an
$\mathbb{F}_{q^n}$-vector space. Its dimension is obviously 1, so
the trace function is an $\mathbb{F}_{q^n}$-basis of
$\mathbb{F}_{q^n}^\vee$ as a nonzero element of it). Define a map
$\psi$ as
\begin{equation}\label{complp}
\begin{split}
\psi:~\mathbb{F}_{q^n}^\vee\otimes_{\mathbb{F}_{q}}\mathbb{F}_{q^n}&\longrightarrow\mathscr{L}_n(\mathbb{F}_{q^n}) \\
T_\alpha\otimes \beta&\longmapsto\tr(\alpha
x)\beta=\sum_{i=0}^{n-1}\beta\alpha^{q^i}x^{q^i}.
\end{split}
\end{equation}
It is easy to see that $\psi$ just defines the isomorphism map we
need.

The isomorphism map $\psi$ indicates that every linearized
polynomial $L(x)\in\mathscr{L}_n(\mathbb{F}_{q^n})$ admits some kind
of ``trace representation" like $L(x)=\sum_{i=0}^{s}\tr (\omega_i
x)\theta_i$ for some $\omega_i$, $\theta_i\in\mathbb{F}_{q^n}$,
$0\leq i\leq s$, since every element of
$\mathbb{F}_{q^n}^\vee\otimes_{\mathbb{F}_{q}}\mathbb{F}_{q^n}$ is a
sum of some single tensors. However, this kind of representation is
not unique for $L(x)$. We will talk about this in detail in Section
\ref{secrep}.


\section{Characterize $\mathscr{L}_n(\mathbb{F}_{q^n})$ via the Dickson matrix algebra and further}\label{secdmatalg}

As mentioned before, we call a matrix of the form
$D_L=\left(a_{j-i}^{q^i}\right)\in\mathscr{M}_n(\mathbb{F}_{q^n})$ a
Dickson matrix or a $\sigma$-circulant matrix associated to the
linearized polynomial
 $L(x)=\sum_{i=0}^{n-1}a_ix^{q^i}\in\mathscr{L}_n(\mathbb{F}_{q^n})$. Denote by
$\mathscr{D}_n(\mathbb{F}_{q^n})$ the set of all $n\times n$ Dickson
matrices over $\mathbb{F}_{q^n}$. It can be verified that
$\mathscr{D}_n(\mathbb{F}_{q^n})$ forms an $\mathbb{F}_{q}$-algebra
under operations in the standard manner. It is clear that
$\mathscr{D}_n(\mathbb{F}_{q^n})$ is an $\mathbb{F}_{q}$-vector
space. To verify that $\mathscr{D}_n(\mathbb{F}_{q^n})$ is a ring,
we need only to verify that $\mathscr{D}_n(\mathbb{F}_{q^n})$ is
closed under multiplication of matrices as associativity always
holds for this multiplication. This can be completed by noting that,
for $a_i$, $b_i\in\mathbb{F}_{q^n}$, $0\leq i\leq n-1$,
\begin{equation}\label{dicmat}
\left(a_{j-i}^{q^i}\right)\left(b_{j-i}^{q^i}\right)=
\left(\sum_{k=0}^{n-1}a_{k-i}^{q^i}b_{j-k}^{q^k}\right)=\left(\sum_{k=0}^{n-1}a_k^{q^i}b_{j-i-k}^{q^{k+i}}\right)
=\left(c_{j-i}^{q^i}\right)
\end{equation}
where $c_i=\sum_{k=0}^{n-1}a_kb_{i-k}^{q^k}$, $0\leq i\leq n-1$.

As an $\mathbb{F}_{q}$-algebra, $\mathscr{D}_n(\mathbb{F}_{q^n})$ is
an $\mathbb{F}_{q}$-subalgebra of $\mathscr{M}_n(\mathbb{F}_{q^n})$.
However, it is not an $\mathbb{F}_{q^n}$-algebra while
$\mathscr{M}_n(\mathbb{F}_{q^n})$ is. Due to its connections to
linearized polynomials, we have the $\mathbb{F}_{q}$-algebra
isomorphism in the following theorem.
\begin{thm}\label{dicalg}
\[\mathscr{L}_n(\mathbb{F}_{q^n})\cong \mathscr{D}_n(\mathbb{F}_{q^n}).\]
\end{thm}
\proof A map from $\mathscr{L}_n(\mathbb{F}_{q^n})$ to
$\mathscr{D}_n(\mathbb{F}_{q^n})$ can be constructed
straightforwardly as
\begin{eqnarray*}
\varphi:~\mathscr{L}_n(\mathbb{F}_{q^n})&\longrightarrow& \mathscr{D}_n(\mathbb{F}_{q^n})\\
L(x)=\sum_{i=0}^{n-1}a_ix^{q^i}&\longmapsto&D_L.
\end{eqnarray*}
$\varphi$ is clearly an isomorphism between vector spaces. On the
other hand, since
\begin{eqnarray*}
  L_1(x)\circ L_2(x) &=& \sum_{i=0}^{n-1}a_i\big(\sum_{j=0}^{n-1}b_jx^{q^j}\big)^{q^i} \\
   &=&\sum_{i=0}^{n-1}\sum_{j=0}^{n-1}a_ib_j^{q^i}x^{q^{i+j}}  \\
   &=&\sum_{i=0}^{n-1}\big(\sum_{k=0}^{n-1}a_kb_{i-k}^{q^k}\big)x^{q^i}
\end{eqnarray*}
for $L_1(x)=\sum_{i=0}^{n-1}a_ix^{q^i}$,
$L_2(x)=\sum_{i=0}^{n-1}b_ix^{q^i}\in\mathscr{L}_n(\mathbb{F}_{q^n})$.
Hence we have
\begin{eqnarray*}
  \varphi(L_1(x)\circ L_2(x)) &=& D_{L_1\circ L_2}=\left(\big(\sum_{k=0}^{n-1}a_kb_{j-i-k}^{q^k}\big)^{q^i}\right) \\
   &=&\left(a_{j-i}^{q^i}\right)\left(b_{j-i}^{q^i}\right)  \\
   &=&\varphi(L_1)\varphi(L_2)
\end{eqnarray*}
form (\ref{dicmat}). Thus $\varphi$ is furthermore an algebra
isomorphism. \qedd

From Theorem \ref{dicalg}, we know that
$\dim_{\mathbb{F}_{q}}\mathscr{D}_n(\mathbb{F}_{q^n})=n^2$. In fact,
we can tell more about the isomorphism in Theorem \ref{dicalg} by
examining the matrix representations of linearized polynomials under
a given basis.

\begin{lem}\label{matrepp}
Let
$L(x)=\sum_{i=0}^{n-1}a_ix^{q^i}\in\mathscr{L}_n(\mathbb{F}_{q^n})$
and $M_L\in\mathscr{M}_n(\mathbb{F}_{q})$ be the matrix of the
linear transformation induced by $L(x)$ under the basis
$\{\beta_i\}_{i=0}^{n-1}$. Then
\begin{equation}\label{matrep}
M_L=\left(\beta_j^{q^i}\right)^{-1}D_L\left(\beta_j^{q^i}\right).
\end{equation}
\end{lem}

\proof Assume $M_L=\big(m_{ij}\big)$. Since
\[\big(L(\beta_0),\ldots,L(\beta_{n-1})\big)=(\beta_0, \ldots, \beta_{n-1})M_L,\]
we have
\[m_{ij}=\tr(\beta_i^*L(\beta_j)), ~~0\leq i,~j\leq n-1,\]
\noindent i.e.
\[M_L=\big(\tr(\beta_i^*L(\beta_j))\big)=\left(\beta_i^{*q^j}\right)\left(L(\beta_j)^{q^i}\right).\]
Note that
\[L(x)^{q^i}=\big(\sum_{j=0}^{n-1}a_jx^{q^j}\big)^{q^i}=\sum_{j=0}^{n-1}a_{j-i}^{q^i}x^{q^j}\]
for $0\leq i\leq n-1$, i.e.
\begin{equation*}\label{linmat}
\begin{pmatrix}
L(x)\\L(x)^q\\\vdots\\L(x)^{q^{n-1}}
\end{pmatrix}=\begin{pmatrix}
a_0&a_1&\dots&a_{n-1}\\
a_{n-1}^q&a_0^q&\dots&a_{n-2}^q\\
\vdots&\vdots&&\vdots\\
a_1^{q^{n-1}}&a_2^{q^{n-1}}&\dots&a_0^{q^{n-1}}
\end{pmatrix}\begin{pmatrix}
x\\x^q\\\vdots\\x^{q^{n-1}}
\end{pmatrix}=D_L\begin{pmatrix}
x\\x^q\\\vdots\\x^{q^{n-1}}
\end{pmatrix}.
\end{equation*}
Hence
\[\left(L(\beta_j)^{q^i}\right)=D_L\left(\beta_j^{q^i}\right).\]
Besides, it is obvious that
$\left(\beta_i^{*q^j}\right)=\left(\beta_j^{q^i}\right)^{-1}$ since
\[\left(\beta_i^{*q^j}\right)\left(\beta_j^{q^i}\right)=\big(\tr(\beta_i^*\beta_j)\big)=I_n,\]
where $I_n$ is the $n\times n$ identity matrix.  At last we get
\begin{equation*}
M_L=\left(\beta_j^{q^i}\right)^{-1}D_L\left(\beta_j^{q^i}\right).
\end{equation*}
\qedd

\begin{prop}\label{gdickson}
For any $L(x)\in\mathscr{L}_n(\mathbb{F}_{q^n})$, $\rk L=\rk D_L$
and $\det L=\det D_L$, where $\rk L$ and $\det L$ are the rank and
determinant, respectively, of the linear transformation induced by
$L(x)$.
\end{prop}
\proof From Lemma \ref{matrepp}, we can obviously get
\[\det L=\det M_L=\det D_L.\]
Besides, we have $\rk M_L=\rk D_L$ though $M_L$ and $D_L$ are
matrices over different fields, since if there exist $P$, $Q\in
GL_n(\mathbb{F}_{q})$, where $GL_n(\mathbb{F}_{q})$ is the general
linear group over $\mathbb{F}_{q}$, such that $M_L=PEQ$, where
\begin{equation}\label{elemat}
E=\begin{array}{c@{\hspace{-5pt}}l}
\begin{pmatrix}
1&&&&&\\&\ddots&&&&\\&&1&&&\\&&&0&&\\&&&&\ddots&\\&&&&&0
\end{pmatrix}&
\begin{array}{l}\left.\rule{0mm}{9mm}\right\}\rk M_L~
\text{rows}\\\rule{0mm}{15mm}
\end{array}
\end{array},
\end{equation}
then
\[D_L=\left(\beta_j^{q^i}\right)PEQ\left(\beta_j^{q^i}\right)^{-1}=RES\]
where $R$, $S\in GL_n(\mathbb{F}_{q^n})$.\qedd

Proposition \ref{gdickson} is a direct generalization of Dickson's
well known criterion on linearized permutation polynomials, which
concludes that a linearized polynomial is a permutation polynomial
if and only if the Dickson matrix associated to it is non-singular
\cite[Chap. 7]{lidl}. It is also implied by Proposition
\ref{gdickson} that the determinant of any Dickson matrix over
$\mathbb{F}_{q^n}$ is an element of $\mathbb{F}_{q}$. A direct
consequence of this when $q=2$ is as follows.

\begin{cor}
Let $L(x)\in\mathscr{L}_n(\mathbb{F}_{2^n})$ be a linearized
permutation polynomial. Then $\det L=1$.
\end{cor}

Recall the isomorphism in Theorem \ref{dicalg} and Theorem
\ref{isomatalg}, we are clear that the set of all non-singular
Dickson matrices form a group, which is isomorphic to the general
linear group $GL_n(\mathbb{F}_{q})$. A direct consequence is that
the inverse of an non-singular Dickson matrix is also of the Dickson
type. For $D_L$ associated to a linearized permutation polynomial
$L(x)$, $D_L^{-1}$ is just the Dickson matrix associated to
$L^{-1}(x)$, the composition inverse of $L(x)$. This supplies us a
method to compute inverse polynomials of linearized permutation
polynomials.

\begin{thm}\label{invpoly}
Let
$L(x)=\sum_{i=0}^{n-1}a_ix^{q^i}\in\mathscr{L}_n(\mathbb{F}_{q^n})$
be a linearized permutation polynomial and $D_L$ be its associated
Dickson matrix. Assume $\bar{a}_i$ is the $(i,0)$-th cofactor of
$D_L$, $0\leq i\leq n-1$. Then $\det
L=\sum_{i=0}^{n-1}a_{n-i}^{q^i}\bar{a}_i$ and
\[L^{-1}(x)=\frac{1}{\det L}\sum_{i=0}^{n-1}\bar{a}_ix^{q^i}=\big(\sum_{i=0}^{n-1}\bar{a}_ix^{q^i}\big)\circ\big(\frac{x}{\det L}\big).\]
\end{thm}

\proof The statement on $\det L$ is obvious. Besides, $\bar{a}_0$,
$\bar{a}_1$, $\ldots$, $\bar{a}_{n-1}$ are entries of the first row
of the adjugate matrix $D_L^\star$ of $D_L$. From linear algebra, we
know that \[D_L^{-1}=\frac{1}{\det L}D_L^\star,\] which is the
Dickson matrix associated to $L^{-1}(x)$. Since $\det
L\in\mathbb{F}_{q}$, the result is straightforward. \qedd

From Theorem \ref{invpoly}, we know that the main task is to compute
cofactors of elements of the first column of a Dickson matrix as its
determinant can be obtained by Laplacian expansion afterwards
 in computing inverse polynomial of a linearized
permutation polynomial. The proof of Theorem \ref{invpoly} also
implies that the adjugate matrix of a non-singular Dickson matrix is
also of the Dickson type. In fact, this holds for any Dickson
matrices.

\begin{lem}\label{dicmatfac}
Let $D\in\mathscr{D}_n(\mathbb{F}_{q^n})$. Then there exist two sets
of elements $\{\alpha'_i\}_{i=0}^{n-1}\subseteq\mathbb{F}_{q^n}$ and
$\{\alpha_i\}_{i=0}^{n-1}\subseteq\mathbb{F}_{q^n}$ such that
\[D=\left(\beta_j^{q^i}\right)\left(\alpha_i^{'q^j}\right)=\left(\alpha_j^{q^i}\right)\left(\beta_i^{q^j}\right).\]
\end{lem}
\proof Assume $D=D_L$ is the Dickson matrix associated to a
linearized polynomial $L(x)\in\mathscr{L}_n(\mathbb{F}_{q^n})$ and
$M_L$ is the matrix of the linear transformation induced by $L(x)$
under the basis $\{\beta_i\}_{i=0}^{n-1}$. From Lemma \ref{matrepp},
we have
\[D=\left(\beta_j^{q^i}\right)M_L\left(\beta_j^{q^i}\right)^{-1}=\left(\beta_j^{q^i}\right)\left(\alpha_i^{'q^j}\right)\]
where
$(\alpha'_0,\ldots,\alpha'_{n-1})=(\beta_0^*,\ldots,\beta_{n-1}^*)M_L^\tau$.
The other part of the lemma can be got in the same way by replacing
$\{\beta_i\}_{i=0}^{n-1}$ by $\{\beta_i^*\}_{i=0}^{n-1}$ in Lemma
\ref{matrepp}.\qedd

\begin{rmk}\label{rankdmat}
From \cite[Lemma 3.51]{lidl}, we know that the determinant of a
Dickson matrix $D$ can be represented in the form
\[\det D=\alpha_0\beta_0\prod_{i=0}^{n-2}
\prod_{c_0,\ldots,c_i\in\mathbb{F}_{q}}\big(\alpha_{i+1}-\sum_{j=1}^ic_j\alpha_j\big)\big(\beta_{i+1}-\sum_{j=1}^ic_j\beta_j\big).\]
Besides, it is obvious from Lemma \ref{dicmatfac} that
\[\rk D=\rk\left(\alpha_j^{'q^i}\right)=\rk\left(\alpha_j^{q^i}\right), \]
where the  notations are the same as in Lemma \ref{dicmatfac}.
\end{rmk}

\begin{lem}\label{mooradj}
Let $\{\alpha_i\}_{i=0}^{n-1}\subseteq\mathbb{F}_{q^n}$.  Let
$\tilde{\alpha}_i$ be the $(0,i)$-th  cofactor of the matrix
$\left(\alpha_j^{q^i}\right)$, $0\leq i\leq n-1$. Then
\[\left(\alpha_j^{q^i}\right)^\star=\left(\tilde{\alpha}_i^{q^j}\right)S,\]
where
$S=\mbox{\rm{diag}}\big((-1)^0,(-1)^{n+1},(-1)^{2(n+1)},\ldots,(-1)^{(n-1)(n+1)}\big),$
i.e. $S=I_n$ when $n$ is odd and
$S=\mbox{\rm{diag}}(1,-1,1,-1,\ldots,1,-1)$ when $n$ is even.
\end{lem}
\proof Let $A=\left(\alpha_j^{q^i}\right)$ and $A_{ij}$ be the
$(i,j)$-th cofactor of $A$, $0\leq i,~j\leq n-1$. We need only to
prove $A_{ji}=(-1)^{j(n+1)}\tilde{\alpha}_i^{q^j}$ for any $0\leq
i,~j\leq n-1$. This can be completed by noting that
\begin{eqnarray*}
  \tilde{\alpha}_i^{q^j} &=& (-1)^{i+2}\bigg(\det \begin{pmatrix}
\alpha_0^q&\cdots&\alpha_{i-1}^q&\alpha_{i+1}^q&\cdots&\alpha_{n-1}^q\\
\alpha_0^{q^2}&\cdots&\alpha_{i-1}^{q^2}&\alpha_{i+1}^{q^2}&\cdots&\alpha_{n-1}^{q^2}\\
\vdots&&\vdots&\vdots&&\vdots\\
\alpha_0^{q^{n-1}}&\cdots&\alpha_{i-1}^{q^{n-1}}&\alpha_{i+1}^{q^{n-1}}&\cdots&\alpha_{n-1}^{q^{n-1}}
\end{pmatrix}\bigg)^{q^j} \\
   &=&(-1)^{i+2}\det \begin{pmatrix}
\alpha_0^{q^{j+1}}&\cdots&\alpha_{i-1}^{q^{j+1}}&\alpha_{i+1}^{q^{j+1}}&\cdots&\alpha_{n-1}^{q^{j+1}}
\\
\vdots&&\vdots&\vdots&&\vdots\\
\alpha_0^{q^{n-1}}&\cdots&\alpha_{i-1}^{q^{n-1}}&\alpha_{i+1}^{q^{n-1}}&\cdots&\alpha_{n-1}^{q^{n-1}}\\
\alpha_0&\cdots&\alpha_{i-1}&\alpha_{i+1}&\cdots&\alpha_{n-1}\\
\vdots&&\vdots&\vdots&&\vdots\\
\alpha_0^{q^{j-1}}&\cdots&\alpha_{i-1}^{q^{j-1}}&\alpha_{i+1}^{q^{j-1}}&\cdots&\alpha_{n-1}^{q^{j-1}}
\end{pmatrix}  \\
   &=&(-1)^{i+2}(-1)^{j(n-j-1)}\det\begin{pmatrix}
   \alpha_0&\cdots&\alpha_{i-1}&\alpha_{i+1}&\cdots&\alpha_{n-1}\\
\vdots&&\vdots&\vdots&&\vdots\\
\alpha_0^{q^{j-1}}&\cdots&\alpha_{i-1}^{q^{j-1}}&\alpha_{i+1}^{q^{j-1}}&\cdots&\alpha_{n-1}^{q^{j-1}}\\
\alpha_0^{q^{j+1}}&\cdots&\alpha_{i-1}^{q^{j+1}}&\alpha_{i+1}^{q^{j+1}}&\cdots&\alpha_{n-1}^{q^{j+1}}
\\
\vdots&&\vdots&\vdots&&\vdots\\
\alpha_0^{q^{n-1}}&\cdots&\alpha_{i-1}^{q^{n-1}}&\alpha_{i+1}^{q^{n-1}}&\cdots&\alpha_{n-1}^{q^{n-1}}
\end{pmatrix}\\
&=&(-1)^{j(n+1)}A_{ji},
\end{eqnarray*}
since $(-1)^{i+2}(-1)^{j(n-j-1)}=(-1)^{j(n+1)+(i+1)+(j+1)}$.\qedd

\begin{cor}\label{cofacbasis}
Notations as in Lemma \ref{mooradj}. If $\{\alpha_i\}_{i=0}^{n-1}$
is a basis of $\mathbb{F}_{q^n}$ over $\mathbb{F}_{q}$, then
$\{\tilde{\alpha}_i\}_{i=0}^{n-1}$ is also a basis.
\end{cor}
\proof From \cite[Lemma 3.51]{lidl}, we know that when
$\{\alpha_i\}_{i=0}^{n-1}$ is a basis,
$\left(\alpha_j^{q^i}\right)$, and consequently
$\left(\alpha_j^{q^i}\right)^\star$, is non-singular. Hence
$\left(\tilde{\alpha}_i^{q^j}\right)$ is non-singular from Lemma
\ref{mooradj}, which implies that $\{\tilde{\alpha}_i\}_{i=0}^{n-1}$
is a basis using \cite[Lemma 3.51]{lidl} again. \qedd

\begin{lem}\label{dicmats}
Let $D\in\mathscr{D}_n(\mathbb{F}_{q^n})$. Then
$SDS\in\mathscr{D}_n(\mathbb{F}_{q^n})$, where $S$ is the matrix
defined in Lemma \ref{mooradj}.
\end{lem}
\proof We only need to consider the case that $n$ is even. Assume
$D$ is a Dickson matrix with entries of the first row
$a_0,~a_1,\ldots,a_{n-1}$. Then it can be verified that $SDS$ is a
Dickson matrix with entries of the first row
$a_0,~-a_1,~a_2,~-a_3,\ldots,a_{n-2},~-a_{n-1}$.\qedd

The following lemma is well known in linear algebra.
\begin{lem}\label{adjprop}
Let $N,~N_1,~N_2\in\mathscr{M}_n(\mathbb{K})$, where $\mathbb{K}$ is
a field. Then $N^{\tau\star}=N^{\star\tau}$,
$(N_1N_2)^\star=N_2^\star N_1^\star$.
\end{lem}

\begin{thm}\label{dicadj}
Let $D\in\mathscr{D}_n(\mathbb{F}_{q^n})$. Then
$D^\star\in\mathscr{D}_n(\mathbb{F}_{q^n})$.
\end{thm}
\proof From Lemma \ref{dicmatfac}, we can assume
$D=\left(\beta_j^{q^i}\right)\left(\alpha_i^{q^j}\right)$ for a set
of elements $\{\alpha_i\}_{i=0}^{n-1}\subseteq\mathbb{F}_{q^n}$.
Then
\begin{eqnarray*}
  D^\star &=& \left(\alpha_i^{q^j}\right)^\star\left(\beta_j^{q^i}\right)^\star \\
   &=&\big(\left(\tilde{\alpha}_i^{q^j}\right)S\big)^\tau\big(\left(\tilde{\beta}_i^{q^j}\right)S\big)\\
   &=&S\left(\tilde{\alpha}_j^{q^i}\right)\left(\tilde{\beta}_i^{q^j}\right)S
\end{eqnarray*}
from Lemma \ref{mooradj} and Lemma \ref{adjprop}, where
$\tilde{\alpha}_i$ and $\tilde{\beta}_i$ are the $(0,i)$-th
cofactors of the matrices $\left(\alpha_j^{q^i}\right)$ and
$\left(\beta_j^{q^i}\right)$ respectively, $0\leq i\leq n-1$. Hence
$D^\star\in\mathscr{D}_n(\mathbb{F}_{q^n})$ from Corollary
\ref{cofacbasis}, Lemma \ref{dicmatfac} and Lemma
\ref{dicmats}.\qedd

Base on Theorem \ref{dicadj}, we introduce the following definition.

\begin{defi}\label{adjpolydef}
Let $L(x)\in\mathscr{L}_n(\mathbb{F}_{q^n})$. The adjugate
polynomial $L^\star(x)$ of $L(x)$ is defined to be the linearized
polynomial associated to $D_L^\star$.
\end{defi}

In fact, the coefficients of $L^\star(x)$ are just cofactors of
elements of the first column of $D_L$. Hence Definition
\ref{adjpolydef} is only meaningful for $L(x)$ with $\rk L\geq n-1$
since $D_L^\star=0$ in other cases. From the fact that $D_L
D_L^\star=D_L^\star D_L=(\det D_L) I_n$, we have
\begin{equation}\label{adjpolycomp}
L(x)\circ L^\star(x)=L^\star(x)\circ L(x)=(\det L)x.
\end{equation}
The case $\rk L=n$ has already been discussed in Theorem
\ref{invpoly} which indicates that $L^{-1}(x)= \DF{1}{\det L}L^\star
(x)$. When $\rk L=n-1$, we have the following property.

\begin{prop}\label{adjpolyprop}
Assume $\rk L=n-1$. Then $\rk L^\star=1$. Furthermore,
\[\ker L^\star=\im L,~\im L^\star=\ker L,\]
where $\ker$ and $\im$ represent kernels and images respectively of
linearized polynomials.
\end{prop}
\proof From linear algebra we know that $\rk D_L^\star=1$ when $\rk
D_L=n-1$ \cite{mortici}, which implies $\rk L^\star=1$. Furthermore,
we have $L(x)\circ L^\star(x)=L^\star(x)\circ L(x)=0$ since $\det
L=0$. Thus $\im L^\star\subseteq\ker L$ and $\im L\subseteq\ker
L^\star$. By comparing dimensions, we know that $\im L^\star=\ker L$
and $\ker L^\star=\im L$.\qedd

Proposition \ref{adjpolyprop} gives characterizations of linearized
polynomials with 1-dimensional kernels. For $L(x)$ with $\rk L=n-1$,
$\ker L=\gamma\cdot\mathbb{F}_{q}$ for some
$\gamma\in\mathbb{F}_{q^n}^*$ and $\im L=\mathcal {H}_\delta$ for
some $\delta\in\mathbb{F}_{q^n}^*$, where $\mathcal
{H}_\delta=\{x\in\mathbb{F}_{q^n}|\;\tr(\delta x)=0\}$ is a
hyperplane in $\mathbb{F}_{q^n}$. $\delta$ can be obtained by
computing $L^\star(x)$.

\begin{eg}\label{ker1}
Let $L(x)=x^q-\gamma^{q-1}x\in\mathscr{L}_n(\mathbb{F}_{q^n})$,
$\gamma\in\mathbb{F}_{q^n}^*$. It is obvious that $\ker
L=\gamma\cdot\mathbb{F}_{q}$. Besides, it is easy to get
\[L^\star(x)=(-1)^{n-1}\gamma\tr\big(\frac{x}{\gamma^q}\big),\]
since the $(i,0)$-th cofactor of the matrix
\[D_L=\left(\begin{array}{cccc}
  -\gamma^{q-1}&1\\
  &-\gamma^{(q-1)q}&1\\
  &&{\ddots}\text{\textcolor{white}{.\;\;\quad.}}&\ddots\text{\textcolor{white}{.\;\;\;\;\quad.}}\\
 &&-\gamma^{(q-1)q^{n-2}}&1\\
  1&&&-\gamma^{(q-1)q^{n-1}}
\end{array}\right)\]
is $(-1)^{n-1}\gamma^{1-q^{i+1}}$, $0\leq i\leq n-1$. Hence $\im
L=\mathcal {H}_{{1}/{\gamma^q}}$.
\end{eg}

Generally, we know from Remark \ref{skewgcrd} that for
$L(x)\in\mathscr{L}_n(\mathbb{F}_{q^n})$ of degree $q^d$, $0\leq
d\leq n-1$, with $\rk L=n-1$ and $\ker
L=\gamma_0\cdot\mathbb{F}_{q}$ for some
$\gamma_0\in\mathbb{F}_{q^n}^*$, there exists
$L_1(x)\in\mathscr{L}_n(\mathbb{F}_{q^n})$ of degree $q^{d-1}$ such
that \[L(x)=L_1(x)\circ\left(x^q-\gamma_0^{q-1}x\right).\] It is
clear that $\rk L_1\geq n-1$ and when $\rk L_1=n-1$,
\[\im (x^q-\gamma_0^{q-1}x)\cap\ker L_1=\{0\}.\]
Assume $\rk L_1=n-1$ and $\ker L_1=\gamma_1\cdot\mathbb{F}_{q}$,
$\gamma_1\in\mathbb{F}_{q^n}^*$. From Example \ref{ker1} we need
$\gamma_1\not\in\mathcal {H}_{{1}/{\gamma_0^q}}$. Inductively, we
finally know that there exist
$\gamma_1,~\gamma_2,~\ldots,~\gamma_{r-1}\in\mathbb{F}_{q^n}^*$ for
some $r\leq d$, satisfying that $\gamma_i\not\in\mathcal
{H}_{{1}/{\gamma_{i-1}^q}}$ (i.e.
$\tr(\frac{\gamma_i}{\gamma_{i-1}^q})\neq0$) for $1\leq i\leq r-1$,
and $L_r(x)\in\mathscr{L}_n(\mathbb{F}_{q^n})$ of degree $q^{d-r}$
which is a linearized permutation polynomial, such that
\[L(x)=L_r(x)\circ\left(x^q-\gamma_{r-1}^{q-1}x\right)\circ\cdots\circ\left(x^q-\gamma_1^{q-1}x\right)
\circ\left(x^q-\gamma_0^{q-1}x\right).\] This can provide an answer
to the open problem proposed in \cite{charp} to an extent. By Lemma
\ref{adjprop} and Example \ref{ker1}, we have
\begin{eqnarray*}
  L^\star(x) &=&\left(x^q-\gamma_0^{q-1}x\right)^\star \circ\left(x^q-\gamma_1^{q-1}x\right)^\star \circ\cdots\circ
  \left(x^q-\gamma_{r-1}^{q-1}x\right)^\star\circ L^\star_r(x)\\[.1cm]
   &=&\left[(-1)^{n-1}\gamma_0\tr\big(\frac{x}{\gamma_0^q}\big)\right]\circ
   \left[(-1)^{n-1}\gamma_1\tr\big(\frac{x}{\gamma_1^q}\big)\right]\circ\cdots\\[.1cm]
   &&\circ
    \left[(-1)^{n-1}\gamma_{r-1}\tr\big(\frac{x}{\gamma_{r-1}^q}\big)\right]\circ
    L^\star_r(x)\\[.1cm]
   &=&(-1)^{r(n-1)}\gamma_0\tr\big(\frac{\gamma_1}{\gamma_0^q}\big)\tr\big(\frac{\gamma_2}{\gamma_1^q}\big)
   \cdots\tr\big(\frac{\gamma_{r-1}}{\gamma_{r-2}^q}\big)\tr\big(\frac{L^\star_r(x)}{\gamma_{r-1}^q}\big).
\end{eqnarray*}


\section{Representations of linearized polynomials}\label{secrep}

Generally we represent a linearized polynomial by giving its
coefficients as a polynomial. Only recently, a new kind of
representation was proposed in \cite{kzhou, pzyuan, lings}.
\begin{thm}[\cite{lings}]\label{replp}
Let $L(x)\in\mathscr{L}_n(\mathbb{F}_{q^n})$ be a linearized
polynomial of rank $k$, where $k$ is an integer, $0\leq k\leq n$.
Then

\noindent (1) there exists a uniquely determined ordered set
$\{\alpha'_i\}_{i=0}^{n-1}\subseteq\mathbb{F}_{q^n}$ of rank $k$
such that $L(x)$ can be represented as
\[L(x)=\sum_{i=0}^{n-1}\tr (\alpha'_i x)\beta_i;\]

\noindent (2) there exists a uniquely determined ordered set
$\{\alpha_i\}_{i=0}^{n-1}\subseteq\mathbb{F}_{q^n}$ of rank $k$ such
that $L(x)$ can be represented as
\[L(x)=\sum_{i=0}^{n-1}\tr (\beta_i x)\alpha_i;\]

\noindent(3) there exist two sets of elements
$\{\omega_i\}_{i=0}^{k-1}\subseteq\mathbb{F}_{q^n}$ and
$\{\theta_i\}_{i=0}^{k-1}\subseteq\mathbb{F}_{q^n}$ both of rank $k$
such that $L(x)$ can be represented as
\[L(x)=\sum_{i=0}^{k-1}\tr (\omega_i
x)\theta_i.\]
\end{thm}

This kind of "trace representations" for linearized polynomials was
firstly motivated by K. Zhou in \cite{kzhou} and made clear by P.Z.
Yuan et al. in \cite{pzyuan} for linearized permutation polynomials.
Theorem \ref{replp} was proved by S. Ling and L.J. Qu in
\cite{lings}, solving an open problem proposed in \cite{charp}
asking for linearized polynomials of rank $n-1$.

In the following, we reestablish Theorem \ref{replp} via the
isomorphism we have constructed in Theorem \ref{compo} and Theorem
\ref{dicalg}. We will find that both approaches are more simple than
that in \cite{lings}. Moreover, we can make the reason for the
existence of such kind of representations for linearized polynomials
clear. Besides, we will talk further about this kind of
representations.

\subsection{Composition algebra approach to Theorem \ref{replp}}\label{sseccomp}
To obtain Theorem \ref{replp} from Theorem \ref{compo}, we firstly
recall the definition of tensor ranks of elements in a tensor space.

\begin{defi}
Let $V=\bigotimes_{i=1}^r V_i$, where $V_i$, $1\leq i\leq r$, are
all vector spaces over a field $\mathbb{K}$. For any $T\in V$, the
tensor rank of it, denoted by $\hbox{\rm{trk}\,}T$, is defined to be
\[\hbox{\rm{trk}\,}T=\min\{k\in\mathbb{N}|\,\exists v_{ji}\in V_j,~1\leq j\leq r,~1\leq i\leq k,
~\text{s.t.}~T=\sum_{i=1}^k v_{1i}\otimes\cdots\otimes v_{ri}\}.\]
\end{defi}

Now considering the tensor space
$\mathbb{F}_{q^n}^\vee\otimes_{\mathbb{F}_{q}}\mathbb{F}_{q^n}$, we
have the following lemma.
\begin{lem}\label{trank}
For any
$T\in\mathbb{F}_{q^n}^\vee\otimes_{\mathbb{F}_{q}}\mathbb{F}_{q^n}$,
\[\hbox{\rm{trk}\,}T=\rk \psi(T),\]
where $\psi$ is the map defined in (\ref{complp}) in Section
\ref{seccompalg}.
\end{lem}
\proof Assume $\hbox{\rm{trk}\,}T=k$ for some positive integer $k$,
then there exist two sets of elements
$\{\omega_i\}_{i=0}^{k-1}\subseteq\mathbb{F}_{q^n}$ and
$\{\gamma_i\}_{i=0}^{k-1}\subseteq\mathbb{F}_{q^n}$ such that
\[T=\sum_{i=0}^{k-1}T_{\omega_i}\otimes \gamma_i.\]
Firstly we prove that
$\rk_{\mathbb{F}_{q}}\{\omega_i\}_{i=0}^{k-1}=\rk_{\mathbb{F}_{q}}\{\gamma_i\}_{i=0}^{k-1}=k$.
If not so, say $\rk_{\mathbb{F}_{q}}\{\gamma_i\}_{i=0}^{k-1}=r<k$,
we have $\gamma_{k-1}=\sum_{i=0}^{k-2}c_i\gamma_i$ for some
$c_i\in\mathbb{F}_{q}$, $0\leq i\leq k-2$, without loss of
generality. Then
\begin{eqnarray*}
   T&=& \sum_{i=0}^{k-2}T_{\omega_i}\otimes\gamma_i+T_{\omega_{k-1}}\otimes \sum_{i=0}^{k-2}c_i\gamma_i\\
   &=& \sum_{i=0}^{k-2}(T_{\omega_i}+c_iT_{\omega_{k-1}}) \otimes\gamma_i\\
   &=&\sum_{i=0}^{k-2}T_{\omega_i+c_i\omega_{k-1}}\otimes\gamma_i,
\end{eqnarray*}
which contradicts the assumption that $\hbox{\rm{trk}\,}T=k$ from
the definition of tensor rank. A similar contradiction can be
derived if $\rk_{\mathbb{F}_{q}}\{\omega_i\}_{i=0}^{k-1}<k$.

Now we extend $\{\omega_i\}_{i=0}^{k-1}$ to be a basis
$\{\omega_i\}_{i=0}^{n-1}$ of $\mathbb{F}_{q^n}$ over
$\mathbb{F}_{q}$ with dual basis $\{\omega_i^*\}_{i=0}^{n-1}$.
Consider the set of elements $\{\psi(T)(\omega_i^*)\}_{i=0}^{n-1}$
in $\mathbb{F}_{q^n}$, which equals
$\{\gamma_i\}_{i=0}^{k-1}\cup\{0\}$ since
\[\psi(T)(\omega_i^*)=\sum_{j=0}^{k-1}\tr(\omega_j\omega_i^*)\gamma_j=\left\{
\begin{array}{cc}
\gamma_i&\text{if}~0\leq i\leq k-1\\
0&\text{if}~k\leq i\leq n-1
\end{array}\right., \]
we get that $\rk
\psi(T)=\rk_{\mathbb{F}_{q}}\{\gamma_i\}_{i=0}^{k-1}=k$.\qedd

From the comments at the end of Section \ref{seccompalg} and Lemma
\ref{trank}, we get Theorem \ref{replp}(3) directly. We remark that
Theorem \ref{replp}(1) and \ref{replp}(2) are both direct
consequences of Theorem \ref{replp}(3). For example, to get
\ref{replp}(2) from \ref{replp}(3), we assume
\[(\omega_0,\ldots,\omega_{k-1})=(\beta_0,\ldots,\beta_{n-1})M,\]
where $M$ is an $n\times k$ matrix over $\mathbb{F}_{q}$ with
$(i,j)$-th entry $m_{ij}$, $0\leq i\leq n-1$, $0\leq j\leq k-1$,
which is of full column rank. Then
\begin{eqnarray*}
   L(x)&=&\sum_{i=0}^{k-1}\tr (\omega_i
x)\theta_i  \\
   &=& \sum_{i=0}^{k-1}\tr\big(\sum_{j=0}^{n-1}m_{ji}\beta_jx\big)\theta_i\\
   &=&\sum_{j=0}^{n-1}\tr(\beta_jx)\big(\sum_{i=0}^{k-1}m_{ji}\theta_i\big)\\
   &=&\sum_{j=0}^{n-1}\tr(\beta_jx)\alpha_j,
\end{eqnarray*}
where
\[(\alpha_0,\ldots,\alpha_{n-1})=(\theta_0,\ldots,\theta_{k-1})M^\tau.\]
It is clear that $\rk_{\mathbb{F}_{q}}\{\alpha_i\}_{i=0}^{n-1}=k$.
The uniqueness of $\{\alpha_i\}_{i=0}^{n-1}$ can also be easily
checked.

In fact, S. Ling and L.J. Qu proved Theorem \ref{replp}(1) and
\ref{replp}(2) firstly and \ref{replp}(3) based on them in
\cite{lings}. Our approach seem more simple to obtain the same
result. Besides, we find that it is just due to Theorem \ref{compo}
and the fact that every element in a tensor space can be represented
as a sum of single tensors, that there exist ``trace
representations" for linearized polynomials.


\subsection{Dickson matrix algebra approach to Theorem \ref{replp}}\label{ssecdicmat}

We firstly propose the following proposition, which is a direct
generalization of Lemma 3.51 in \cite{lidl}.

\begin{prop}\label{mooredetrk}
Let $k\in\mathbb{N}$ and
$\{\alpha_0,\ldots,\alpha_{k-1}\}\subseteq\mathbb{F}_{q^n}$. Then
\[\rk_{\mathbb{F}_{q}}\{\alpha_0,\ldots,\alpha_{k-1}\}=\rk \Delta(\alpha_0,\ldots,\alpha_{k-1}),\]
where
$\Delta(\alpha_0,\ldots,\alpha_{k-1})=\left(\alpha_j^{q^i}\right)_{0\leq
i\leq k-1,\;0\leq j\leq k-1}$.
\end{prop}
\proof It is clear that
$\rk_{\mathbb{F}_{q}}\{\alpha_0,\ldots,\alpha_{k-1}\}\geq\rk
\Delta(\alpha_0,\ldots,\alpha_{k-1})$. Assume
$\rk_{\mathbb{F}_{q}}\{\alpha_0,\ldots,\alpha_{k-1}\}=r>\rk
\Delta(\alpha_0,\ldots,\alpha_{k-1})=r-1$, and
\[\begin{pmatrix}
\alpha_0&\alpha_1&\dots&\alpha_{r-2}\\
\alpha_0^q&\alpha_1^q&\dots&\alpha_{r-2}^q\\
\vdots&\vdots&&\vdots\\
\alpha_0^{q^{k-1}}&\alpha_1^{q^{k-1}}&\dots&\alpha_{r-2}^{q^{k-1}}
\end{pmatrix}\]
is a submatrix of $\Delta(\alpha_0,\ldots,\alpha_{k-1})$ with full
column rank, without loss of generality. Then for any
$(c_0,\ldots,c_{r-2})\in\mathbb{F}_{q}^{r-1}\backslash \{\vec{0}\}$,
\[\sum_{i=0}^{r-2}c_i\begin{pmatrix}
\alpha_i\\
\alpha_i^q\\
\vdots\\
\alpha_i^{q^{k-1}}
\end{pmatrix}=\begin{pmatrix}
\sum_{i=0}^{r-2}c_i\alpha_i\\
\big(\sum_{i=0}^{r-2}c_i\alpha_i\big)^q\\
\vdots\\
~~\big(\sum_{i=0}^{r-2}c_i\alpha_i\big)^{q^{k-1}}
\end{pmatrix}\neq\begin{pmatrix}
0\\
0\\
\vdots\\
0
\end{pmatrix},\]
thus $\sum_{i=0}^{r-2}c_i\alpha_i\neq 0$, which implies that
$\{\alpha_0,\ldots,\alpha_{r-2}\}$ is  linearly independent over
$\mathbb{F}_{q}$. Since $\rk_{\mathbb{F}_{q}}
\{\alpha_0,\ldots,\alpha_{k-1}\}=r$, we assume
$\{\alpha_0,\ldots,\alpha_{r-2},\alpha_{r-1}\}$ is linear
independent over $\mathbb{F}_{q}$ without loss of generality. Then
the submatrix
\[\begin{pmatrix}
\alpha_0&\alpha_1&\dots&\alpha_{r-1}\\
\alpha_0^q&\alpha_1^q&\dots&\alpha_{r-1}^q\\
\vdots&\vdots&&\vdots\\
\alpha_0^{q^{k-1}}&\alpha_1^{q^{k-1}}&\dots&\alpha_{r-1}^{q^{k-1}}
\end{pmatrix}\]
of $\Delta(\alpha_0,\ldots,\alpha_{k-1})$ is with column rank $r-1$.
Hence there exists
$(c_0,\ldots,c_{r-1})\in\mathbb{F}_{q^n}^{r-1}\backslash
\{\vec{0}\}$ such that
\[\sum_{i=0}^{r-1}c_i\alpha_i^{q^j}=0\]
for any $0\leq j\leq k-1$. We can assume $c_{r-1}=1$ as it is
nonzero. Furthermore, since $\alpha_{r-1}\neq 0$, there exists a
$t$, $0\leq t\leq r-2$, such that $c_t\neq 0$. Now for any $0\leq
j\leq k-2$, we have
\begin{eqnarray*}
 \big(\sum_{i=0}^{r-1}c_i\alpha_i^{q^j}\big)^q-\sum_{i=0}^{r-1}c_i\alpha_i^{q^{j+1}} &=& \sum_{i=0}^{r-1}(c_i^q-c_i)\alpha_i^{q^{j+1}} \\
   &=&\sum_{i=0}^{r-2}(c_i^q-c_i)\alpha_i^{q^{j+1}}.
\end{eqnarray*}
Since the submatrix
\[\begin{pmatrix}
\alpha_0^q&\alpha_1^q&\dots&\alpha_{r-2}^q\\
\alpha_0^{q^2}&\alpha_1^{q^2}&\dots&\alpha_{r-2}^{q^2}\\
\vdots&\vdots&&\vdots\\
\alpha_0^{q^{k-1}}&\alpha_1^{q^{k-1}}&\dots&\alpha_{r-2}^{q^{k-1}}
\end{pmatrix}\]
is with full column rank, we have
\[c_i^q-c_i=0\]
for any $0\leq i\leq r-2$. Hence $c_i\in\mathbb{F}_{q}$ for any
$0\leq i\leq r-1$. However, $\sum_{i=0}^{r-1}c_i\alpha_i=0$
contradicts the fact that $\{\alpha_0,\ldots,\alpha_{r-1}\}$ is
linear independent over $\mathbb{F}_{q}$.\qedd

For an linearized polynomial
$L(x)=\sum_{i=0}^{n-1}a_ix^{q^i}\in\mathscr{L}_n(\mathbb{F}_{q^n})$
with $\rk L=k$, $0\leq k\leq n$,
$D_L=\left(\beta_j^{q^i}\right)\left(\alpha_i^{'q^j}\right)=\left(\alpha_j^{q^i}\right)\left(\beta_i^{q^j}\right)$
for two sets of elements
$\{\alpha'_i\}_{i=0}^{n-1}\subseteq\mathbb{F}_{q^n}$ and
$\{\alpha_i\}_{i=0}^{n-1}\subseteq\mathbb{F}_{q^n}$ from Lemma
\ref{dicmatfac}, which implies that
\begin{eqnarray*}
 L(x)&=&(a_0,a_1,\ldots,a_{n-1})\begin{pmatrix}
   x\\x^q\\\vdots\\x^{q^{n-1}}
 \end{pmatrix}  \\
 &=&(\beta_0,\beta_1,\ldots,\beta_{n-1})\begin{pmatrix}
\alpha'_0&\alpha_0^{'q}&\dots&\alpha_0^{'q^{n-1}}\\
\alpha'_1&\alpha_1^{'q}&\dots&\alpha_1^{'q^{n-1}}\\
\vdots&\vdots&&\vdots\\
\alpha'_{n-1}&\alpha_{n-1}^{'q}&\dots&\alpha_{n-1}^{'q^{n-1}}
\end{pmatrix}\begin{pmatrix}
   x\\x^q\\\vdots\\x^{q^{n-1}}
 \end{pmatrix}\\
 &&\bigg(~=(\alpha_0,\alpha_1,\ldots,\alpha_{n-1})\begin{pmatrix}
\beta_0&\beta_0^q&\dots&\beta_0^{q^{n-1}}\\
\beta_1&\beta_1^q&\cdots&\beta_1^{q^{n-1}}\\
\vdots&\vdots&&\vdots\\
\beta_{n-1}&\beta_{n-1}^q&\dots&\beta_{n-1}^{q^{n-1}}
\end{pmatrix}\begin{pmatrix}
   x\\x^q\\\vdots\\x^{q^{n-1}}
 \end{pmatrix}~\text{respectively}~\bigg) \\
 &=&\sum_{i=0}^{n-1}\tr (\alpha'_i
 x)\beta_i\qquad\bigg(~=\sum_{i=0}^{n-1}\tr (\beta_i x)\alpha_i~\text{respectively}~\bigg).
\end{eqnarray*}
Furthermore, since
$\rk\left(\alpha_j^{'q^i}\right)=\rk\left(\alpha_j^{q^i}\right)=\rk
D_L=k$ from Remark \ref{rankdmat} and Proposition \ref{gdickson}, it
is straightforward that
$\rk_{\mathbb{F}_{q}}\{\alpha'_i\}_{i=0}^{n-1}=\rk_{\mathbb{F}_{q}}\{\alpha_i\}_{i=0}^{n-1}=k$
from Proposition \ref{mooredetrk}. Hence Theorem \ref{replp}(1) and
\ref{replp}(2) are obtained.

Particularly, when $k=1$,
$\rk_{\mathbb{F}_{q}}\{\alpha_i\}_{i=0}^{n-1}=1$, so there exist
$c_0$, $c_1$, $\cdots$, $c_{n-1}\in\mathbb{F}_{q}$ which are not all
zero such that $\{\alpha_i\}_{i=0}^{n-1}=\{c_0\theta, c_1\theta,
\ldots, c_{n-1}\theta\}$ for some $\theta\in\mathbb{F}_{q^n}^*$.
Hence
\[\begin{pmatrix}
\alpha_0&\alpha_1&\dots&\alpha_{n-1}\\
\alpha_0^q&\alpha_1^q&\dots&\alpha_{n-1}^q\\
\vdots&\vdots&&\vdots\\
\alpha_0^{q^{n-1}}&\alpha_1^{q^{n-1}}&\dots&\alpha_{n-1}^{q^{n-1}}
\end{pmatrix}=\begin{pmatrix}
   \theta\\\theta^q\\\vdots\\\theta^{q^{n-1}}
 \end{pmatrix}(c_0,c_1,\ldots,c_{n-1})\]
and furthermore,
\begin{eqnarray*}
D_L&=&\begin{pmatrix}
   \theta\\\theta^q\\\vdots\\\theta^{q^{n-1}}
 \end{pmatrix}(c_0,c_1,\ldots,c_{n-1})\begin{pmatrix}
\beta_0&\beta_0^q&\dots&\beta_0^{q^{n-1}}\\
\beta_1&\beta_1^q&\cdots&\beta_1^{q^{n-1}}\\
\vdots&\vdots&&\vdots\\
\beta_{n-1}&\beta_{n-1}^q&\dots&\beta_{n-1}^{q^{n-1}}
\end{pmatrix}\\
&=&\begin{pmatrix}
   \theta\\\theta^q\\\vdots\\\theta^{q^{n-1}}
 \end{pmatrix}(\omega,\omega^q,\ldots,\omega^{q^{n-1}})
\end{eqnarray*}
where $\omega=\sum_{i=0}^{n-1}c_i\beta_i$, i.e.
$L(x)=\theta\tr(\omega x)$.

Now for any $k$, it is an easy exercise in linear algebra that every
matrix of rank $k$ over a field can be factorized into a sum of $k$
matrices each of which is of rank $1$ (the factorization is not
necessarily unique). Hence for $D_L$ with rank $k$, $0\leq k\leq
n-1$, there exist two sets of elements
$\{\omega_i\}_{i=0}^{k-1}\subseteq\mathbb{F}_{q^n}$ and
$\{\theta_i\}_{i=0}^{k-1}\subseteq\mathbb{F}_{q^n}$ such that
\[D_L=\sum_{i=0}^{k-1}\begin{pmatrix}
   \theta_i\\\theta_i^q\\\vdots\\\theta_i^{q^{n-1}}
 \end{pmatrix}(\omega_i,\omega_i^q,\ldots,\omega_i^{q^{n-1}}),\]
which implies that $L(x)=\sum_{i=0}^{k-1}\tr(\omega_ix)\theta_i$.
$\rk\{\omega_i\}_{i=0}^{k-1}=\rk\{\theta_i\}_{i=0}^{k-1}=k$ can be
easily derived, since if not so, $L(x)$ can be represented as a sum
of less than $k$ trace terms, and hence $D_L$ can be factorized into
a sum of less than $k$ rank 1 matrices, which contradicts $\rk
D_L=k$ due to  the fact that the rank of a sum of some matrices is
no bigger than the sum of their ranks. Thus Theorem \ref{replp}(3)
is obtained.


\subsection{Further remarks on Theorem \ref{replp}}\label{ssecrep}
Now we focus on the representation of a linearized polynomial in
Theorem \ref{replp}(2). For linearized permutation polynomials, the
following propositions can be easily verified.

\begin{prop}\label{linpp}
(1) Assume $\{\alpha_i\}_{i=0}^{n-1}$ is a basis of
$\mathbb{F}_{q^n}$ over $\mathbb{F}_{q}$ with dual basis
$\{\alpha^*_i\}_{i=0}^{n-1}$, and $L(x)=\sum_{i=0}^{n-1}\tr (\beta_i
x)\alpha_i$. Then
\[L^{-1}(x)=\sum_{i=0}^{n-1}\tr (\alpha^*_i x)\beta^*_i;\]
(2)\[x=\sum_{i=0}^{n-1}\tr (\beta_i x)\beta^*_i.\]
\end{prop}

Proposition \ref{linpp}(2) is meaningful in understanding some
concepts. For example, we always call the function $\tr(af(x))$,
$a\in\mathbb{F}_{q^n}$, a component function of a polynomial $f(x)$
over $\mathbb{F}_{q^n}$ in cryptography. From Proposition
\ref{linpp}(2), we have
\[f(x)=\sum_{i=0}^{n-1}\tr (\beta_i^*
f(x))\beta_i.\] The functions $\tr(\beta_i^* f(x))$,
$i=0,\ldots,n-1$, are just the component functions of the map
induced by $f(x)$ with respect to the basis
$\{\beta_i\}_{i=0}^{n-1}$. Obviously, the so-called component
functions of $f(x)$, $\tr(af(x))$, $a\in\mathbb{F}_{q^n}$, are just
all possible $\mathbb{F}_{q}$-linear combinations of the component
functions $\{\tr(\beta_i^* f(x))\}_{i=0}^{n-1}$.

Let $L(x)=\sum_{i=0}^{n-1}\tr (\beta_i
x)\alpha_i\in\mathscr{L}_n(\mathbb{F}_{q^n})$ be a linearized
polynomial of rank $k$. If we assume
\[(\alpha_0,\ldots,\alpha_{n-1})=(\beta_0,\ldots,\beta_{n-1})A_L=(\beta_0^*,\ldots,\beta_{n-1}^*)B_L,\]
where $A_L$, $B_L\in\mathscr{M}_n(\mathbb{F}_{q})$ are both of rank
$k$. Then we can represent $L(x)$ as
\begin{equation}\label{repmat1}
L(x)=(\beta_0,\ldots,\beta_{n-1})A_L\begin{pmatrix}\tr(\beta_0x)\\\vdots\\\tr(\beta_{n-1}x)\end{pmatrix}
\end{equation}
or
\begin{equation}\label{repmat2}
L(x)=(\beta_0^*,\ldots,\beta_{n-1}^*)B_L\begin{pmatrix}\tr(\beta_0x)\\\vdots\\\tr(\beta_{n-1}x)\end{pmatrix}.
\end{equation}
These two kinds of representations for linearized polynomials are
convenient when dealing with some problems. Note that
\[D_L=\big(\beta_j^{q^i}\big)A_L\big(\beta_i^{q^j}\big)=\big(\beta_j^{*q^i}\big)B_L\big(\beta_i^{q^j}\big).\]
It is easy to see from Lemma \ref{matrepp} that $B_L$ is just the
matrix of the linear transformation induced by $L(x)$ under the
basis $\{\beta_i^*\}_{i=0}^{n-1}$. Hence we can easily correspond a
given matrix in $\mathscr{M}_n(\mathbb{F}_{q})$ to a linearized
polynomial in $\mathscr{L}_n(\mathbb{F}_{q^n})$ represented in the
form (\ref{repmat2}). Under this correspondence, the three kinds of
elementary matrices admit three kinds of ``elementary linearized
polynomials" of the following form:\\
(1)\[L_{ij}(x)=x-(\tr(\beta_ix)-\tr(\beta_jx))(\beta_i^*-\beta_j^*);\]
(2)\[L_{a,i}(x)=x+a\tr(\beta_ix)\beta_i^*,~~a\in\mathbb{F}_{q^n};\]
(3)\[L_{i+j}(x)=x+\tr(\beta_ix)\beta_j^*.\] From linear algebra, we
know that every linearized polynomial of rank $k$ can be represented
as compositions of some elementary linearized polynomials and a
linearized polynomial of the form $\sum_{i=0}^{k-1}\tr (\beta_i
x)\beta_i^*$, which is the linearized polynomial corresponds to the
matrix in the form (\ref{elemat}).


\section{Characterizations of $\mathscr
{L}_n(\mathbb{F}_{q^m})$}\label{secsubalg}

The isomorphisms in the following theorem to characterize the
subalgebra $\mathscr {L}_n(\mathbb{F}_{q^m})$ are easily to obtain.
\begin{thm}[\cite{mcdonald}]
\[\mathscr{L}_n(\mathbb{F}_{q^m})\cong
\mathbb{F}_{q^m}[x;\sigma]/(x^n-1)\cong\mathbb{F}_{q^m}[G].\]
\end{thm}

Now we study the matrix of the linear transformation induced by an
element of $\mathscr {L}_n(\mathbb{F}_{q^m})$. We use the notations
in Section \ref{secrep} and assume that $\{\beta_i\}_{i=0}^{n-1}$ is
a normal basis generated by $\beta$ with dual basis generator
$\beta^*$ here. For $L(x)=\sum_{i=0}^{n-1}\tr (\beta^{q^i}
x)\alpha_i\in\mathscr {L}_n(\mathbb{F}_{q^m})$, we have
\[\alpha_i=L(\beta^{*q^i}),~~0\leq i\leq n-1.\]
Assume $i=jm+k$, $j\geq 0$, $0\leq k\leq m-1$. Then
\begin{equation}\label{subalgset}
\alpha_i=L(\beta^{*q^{jm+k}})=\big(L(\beta^{*q^k})\big)^{q^{jm}}=\alpha_k^{q^{jm}}.
\end{equation}
Hence the ordered set $\{\alpha_i\}_{i=0}^{n-1}$ is of the form
\begin{equation}\label{subalgrep}
\{\alpha_0,\ldots,\alpha_{m-1},\alpha_0^{q^m},\ldots,\alpha_{m-1}^{q^m},\ldots,\alpha_0^{q^{(t-1)m}},\ldots,\alpha_{n-1}^{q^{(t-1)m}}\}
\end{equation}
where $n=mt$. Then we have the following theorem.
\begin{thm}
Let $V_m$ be the subalgebra of
$\mathbb{F}_{q^n}^\vee\otimes_{\mathbb{F}_{q}}\mathbb{F}_{q^n}$
formed by elements of the form
\[\sum_{i=0}^{n-1}T_{\beta^{q^i}}\otimes \alpha_i\]
where $\beta$ is a normal basis generator and
$\{\alpha_i\}_{i=0}^{n-1}$ is of the form (\ref{subalgrep}). Then
\[\mathscr
{L}_n(\mathbb{F}_{q^m})\cong V_m.\]
\end{thm}

From (\ref{subalgset}), it can be checked that the matrix $B_L$ in
(\ref{repmat2}) is of the form
\[\begin{pmatrix}
B_0&B_1&\cdots&B_{t-1}\\
B_{t-1}&B_0&\cdots&B_{t-2}\\
\vdots&\vdots&&\vdots\\
B_1&B_2&\cdots&B_0
\end{pmatrix},\]
where $B_i\in\mathscr{M}_m(\mathbb{F}_q)$, $0\leq i\leq t-1$, which
is a block circulant matrix of type $(t,m)$ \cite{davis}. Denote by
$\mathscr{B}\mathscr{C}_{t,m}(\mathbb{F}_q)$ the algebra formed by
all block circulant matrices of  type $(t,m)$ over $\mathbb{F}_q$
which is a subalgebra of $\mathscr{M}_n(\mathbb{F}_q)$. The
following theorem is straightforward form Theorem \ref{isomatalg}.
\begin{thm}
\[\mathscr{L}_n(\mathbb{F}_{q^m})\cong\mathscr{B}\mathscr{C}_{t,m}(\mathbb{F}_q).\]
\end{thm}

Particularly, when $m=1$, $\{\alpha_i\}_{i=0}^{n-1}$ is of the form
$\{\alpha^{q^i}\}_{i=0}^{n-1}$ for some $\alpha\in\mathbb{F}_{q^n}$,
and $B_L$ is an $n\times n$ circulant matrix in this case.

Since there is a natural isomorphism between
$\mathscr{B}\mathscr{C}_{t,m}(\mathbb{F}_q)$ and
$\mathscr{M}_m(\mathbb{F}_q)[x]/(x^n-1)$ \cite{davis}, and
\begin{eqnarray*}
  \mathscr{M}_m(\mathbb{F}_q)[x]/(x^t-1) &\cong&\mathscr{M}_m(\mathbb{F}_q)\otimes_{\mathbb{F}_q}\mathbb{F}_q[x]/(x^t-1) \\
   &\cong&\mathscr{M}_m\big(\mathbb{F}_q[x]/(x^t-1)\big),
\end{eqnarray*}
we rediscover the following result firstly obtained by Brawley et
al. in \cite{brawley}, through a pure matrix theoretic approach.
\begin{cor}
\[\mathscr{L}_n(\mathbb{F}_{q^m})\cong\mathscr{M}_m\big(\mathbb{F}_q[x]/(x^t-1)\big).\]
\end{cor}

It can also be checked that the matrix $D_L$ associated to $L(x)$ is
also a  block circulant matrix of type $(t,m)$. Each block of it is
over $\mathbb{F}_{q^m}$, but the entries of the blocks cannot be any
elements of $\mathbb{F}_{q^m}$ as there are relations between them.
All the Dickson matrices associated to elements in
$\mathscr{L}_n(\mathbb{F}_{q^m})$ form a $\mathbb{F}_{q}$-subalgebra
of $\mathscr{D}_n(\mathbb{F}_{q^n})$, which is isomorphic to
$\mathscr{L}_n(\mathbb{F}_{q^m})$.

Let $L(x)=\sum_{i=0}^{n-1}a_ix^{q^i}=\sum_{i=0}^{n-1}\tr
(\beta^{q^i} x)\alpha_i$, then
\begin{eqnarray*}
 a_i&=& \sum_{l=0}^{n-1}\alpha_l\beta^{q^{i+l}} \\
&=& \sum_{j=0}^{t-1}\big(\sum_{k=0}^{m-1}\alpha_k\beta^{q^{i+k}}\big)^{q^{jm}} \\
 &=&\tr_{n/m}\big(\sum_{k=0}^{m-1}\alpha_k\beta^{q^{i+k}}\big),
\end{eqnarray*}
where $\tr_{n/m}$ is the trace function from $\mathbb{F}_{q^n}$ to
$\mathbb{F}_{q^m}$. Particularly, $L(x)$ can be represented as
\[L(x)=\sum_{i=0}^{n-1}\tr(\alpha\beta^{q^i})x^{q^i}=\sum_{i=0}^{n-1}\tr(\alpha^{q^{n-i}}\beta)x^{q^i}\]
for some $\alpha\in\mathbb{F}_{q^n}$ when $m=1$. We can get the
following interesting property.
\begin{prop}
\[\rk_{\mathbb{F}_{q}}\{\alpha^{q^i}\}_{i=0}^{n-1}=n-\deg \hbox{\rm{gcd}}\big(\sum_{i=0}^{n-1}\tr(\alpha\beta^{q^i})x^{q^i},~x^n-1\big),\]
where $\gcd$ of two polynomials represents their greatest common
divisor.
\end{prop}


\section{Concluding remarks}\label{secconc}

In this paper, we propose two new characterizations of the
linearized polynomial algebra $\mathscr{L}_n(\mathbb{F}_{q^n})$. The
new characterizations can help us to get more results about
linearized polynomials over finite fields. As an example, we can
explain the existence of some special kind of ``trace
representations" of linearized polynomials proposed recently and
rediscover it in more simple ways. We remark that the Dickson matirx
associated to a linearized polynomial is important in studying some
problems, say in \cite{ddlin}, it is used to classify equivalent
classes of polynomials used in MPKC. As a result, the isomorphism we
construct in Section \ref{secdmatalg} is meaningful since we can
transform many problems related to linearized polynomials into
matrix theoretic formats. In fact, we have used it to study
quadratic exponential sums recently, but that is not a topic of this
paper.

\end{document}